\documentclass[12pt]{amsart}

\usepackage{amssymb}
\usepackage{graphicx}
\usepackage{epsf}
\usepackage{amsmath}
\usepackage[latin1]{inputenc}
\usepackage{amsfonts}
\usepackage{color}

\numberwithin{equation}{section}
\newtheorem{theo}{Theorem}[section]
\newtheorem{prop}[theo]{Proposition}

\newtheorem{lemm}[theo]{Lemma}

\newtheorem{rema}[theo]{\it Remark}
\newtheorem{exem}[theo]{\it Example}

\def\X{{\bf x}}

\def\P{{\bf p}}

\def\RR{{\mathcal R}}

\def\M{\hbox{{\bf m}}}

\def\join{|}
\def\sep{.} 

\def\kk{\Bbbk}
\def\one{{\bf 1}}
\def\zero{{\bf 0}}
\def\walg{{(\kk \Pi_n, \wedge)}}
\def\valg{{(\kk \Pi_n, \vee)}}
\def\bPiw{({\mathbf \Pi}, \wedge)}
\def\bPiv{({\bf \Pi}, \vee)}
\def\oper{{@}}

\def\shuffle{{\,\raise
1pt\hbox{$\scriptscriptstyle\cup{\mskip-4mu}\cup$}\,}}

\def\NCSym{{\mathsf{NCSym}}}
\def\Res{{\textnormal{Res}}}
\def\Ind{{\textnormal{Ind}}}

\title[Grothendieck bialgebras, Partition lattices, and $NCSym$ ]{Grothendieck bialgebras, Partition lattices and
symmetric functions in noncommutative variables}

\author{N. Bergeron}\address[Nantel Bergeron]{Department of Mathematics and Statistics\\ York  University\\       To\-ron\-to, Ontario M3J 1P3\\ CANADA} \email{bergeron@mathstat.yorku.ca}  \urladdr{http://www.math.yorku.ca/bergeron}

\author{C. Hohlweg}\address[Christophe Hohlweg]{The Fields Institute\\ 222 College Street\\
Toronto, Ontario, M5T 3J1\\ CANADA} \email{chohlweg@fields.utoronto.ca}
\urladdr{http://www.fields.utoronto.ca/\~{}chohlweg}

 \author{M. Rosas}\address[Mercedes Rosas]{Department of Mathematics and Statistics\\ York  University\\       To\-ron\-to, Ontario M3J 1P3\\ CANADA} \email{rosas@mathstat.yorku.ca}

\author{M. Zabrocki} \address[Mike Zabrocki]{Department of Mathematics and Statistics\\ York  University\\       To\-ron\-to, Ontario M3J 1P3\\ CANADA} \email{zabrocki@mathstat.yorku.ca}  \urladdr{http://www.math.yorku.ca/\~{}zabrocki}

\date{\today}

 \thanks{This work is supported in part by CRC and NSERC. It is the results of a working seminar at Fields Institute with the active participation of T. MacHenry, M. Mishna, H. Li and L. Sabourin. }

\begin{document}

\begin{abstract}
We show that the Grothendieck bialgebra of the semi-tower of partition lattice algebras is isomorphic to the graded dual of the bialgebra of symmetric functions in noncommutative variables. In particular this isomorphism singles out a canonical new basis of the symmetric functions in noncommutative variables which would be an analogue of the Schur function basis for this bialgebra.
\end{abstract}

\maketitle
 \section*{Introduction}

Combinatorial Hopf algebras are graded connected Hopf algebras
equipped with a multiplicative linear functional $\zeta :
{\mathcal H} \rightarrow \kk$ called a character (see \cite{ABS}).
Here we assume that $\kk$ is a field of characteristic zero. There
has been renewed interest in these spaces in recent papers (see
for example \cite{AS, BHoh,BMSW, HivThib, Hoff}  and the
references therein) and one particularly interesting aspect of
recent work has been to realize a given combinatorial Hopf algebra
as the Grothendieck Hopf algebra of a tower of algebras.

The prototypical example is
the Hopf algebra of symmetric functions viewed, via the
Frobenius characteristic map, as the Grothendieck Hopf algebras
of the modules of all symmetric group algebras $\kk S_n$ for
$n\ge 0$. The multiplication is
given via induction from $\kk S_n\otimes \kk S_m$ to $\kk S_{n+m}$
and the comultiplication is the sum over $r$ of the restriction
from $\kk S_n$ to $\kk S_{r}\otimes \kk S_{n-r}$. The tensor product
of modules defines a third operation on symmetric
functions usually referred to as the internal multiplication
or the Kronecker product \cite{mac, Ze}. The
Schur symmetric functions are then canonically defined as
the Frobenius image of the simple modules.

There are many more examples of this kind of connection
(see \cite{BHT,Hiv,NCSF4}). Here we are interested in the
bialgebra structure of the symmetric functions in {\sl noncommutative}
variables \cite{BRRZ,BZ,BeCo,RS,W} and the goal of this
paper is to realize it as the
Grothendieck bialgebra of the modules of the
partition lattice algebras.

We denote by $\NCSym=\bigoplus_{d\ge 0} \NCSym_d$ the algebra of
symmetric functions in noncommutative variables, the product is induced
from  the concatenation of words.  This
is a Hopf algebra equipped with an internal comultiplication.
The space $\NCSym_d$ is the subspace of series in the noncommutative
variables $x_1,x_2, \ldots$ with homogeneous degree $d$
that are invariants by any finite permutation of the variables.
The algebra structure of $\NCSym$ was first introduced in \cite{W}
were it was shown to be a free noncommutative algebra.
This algebra was used in \cite{BeCo} to study free powers of
noncommutative rings.  More recently,  a series of
new bases was given  for this space, lifting some of the classical bases
of (commutative) symmetric functions \cite{RS}. The Hopf algebra structure
was uncovered in \cite{AM, BRRZ,BZ} along with other fundamental
algebraic and geometric structures.

The (external) comultiplication $\Delta\colon \NCSym_d \to
\bigoplus \NCSym_k\otimes \NCSym_{d-k}$ is graded and gives rise
to a structure of a graded Hopf algebra on $\NCSym$. The algebra
$\NCSym$ also has an internal comultiplication
$\Delta^{\odot}\colon \NCSym_d \to \NCSym_d\otimes \NCSym_d$ which
is not graded. The algebra $\NCSym$ with the comultiplication
$\Delta^{\odot}$ is only a bialgebra (not graded) and is different
from the previous graded Hopf structure.

After investigating the Hopf algebra structure of $\NCSym$, it is
natural to ask if there exists a tower of algebras
$\{A_n\}_{n \geq 0}$ such that the Hopf algebra $\NCSym$ corresponds to
the Grothendieck bialgebra (or Hopf) algebra of the $A_n$-modules.
 This was the 2004-2005 question for  our algebraic combinatorics working
 seminar at Fields Institute where the research for this article was done.

Our answer involves the partition lattice algebras $(\kk \Pi_n,\wedge)$ and $(\kk\Pi_n,\vee)$
(as well as the Solomon-Tits algebras \cite{Brown,Schocker,Tits}).
 For each one, with finite modules we can define a tensor product of
$\kk \Pi_n$ modules and a restriction from $\kk \Pi_n$
module to $\kk \Pi_k\otimes \kk \Pi_{n-k}$ modules.
This allows us to place on $\bigoplus_n G_0(\kk \Pi_n)$, the
Grothendieck ring of the $\kk \Pi_n$, a bialgebra structure (but
not a Hopf algebra structure). We then define a bialgebra isomorphism
$\bigoplus_n G_0(\kk \Pi_n)\to \NCSym^*$. We call this map the
Frobenius characteristic map of the partition lattice algebras.
This singles out a unique canonical basis of $\NCSym$ (up to
automorphism) corresponding to the simple modules of
the $\kk \Pi_n$.

Our paper is divided into 4 sections as follows. In section
\ref{sec:basis} we recall the definition and structure of
$\NCSym$. We then state our first theorem claiming the existence
of a basis $\mathbf x$ of $\NCSym$ defined by certain algebraic
properties. The proof of it will be postponed to section
\ref{sec:frob}. In section \ref{sec:groth} we recall the
definition and structure of the partition lattice algebras $\kk
\Pi_n$ with the product given by the lattice operation $\wedge$
and define their modules.  We then introduce a structure of a
semi-tower of algebras (i.e. we have a non-unital embedding
$\rho_{n,m} : \kk \Pi_n\otimes \kk\Pi_m \to \kk\Pi_{n+m}$ of
algebras) on the partition lattice algebras and show that it
induces a bialgebra structure on its Grothendieck ring. Our second
theorem states that this Grothendieck bialgebra is dual to
$\NCSym$.  The classes of simple modules correspond then to the
basis $\mathbf x$.
  In view of the work of
Brown \cite{Brown} we remark that this can also be done with the
semi-tower of Solomon-Tits algebras. In section \ref{sec:grothII}
we build the same construction with the lattice algebras $\kk
\Pi_n$ with the product $\vee$.  With this tower  of algebras
(i.e. $\rho_{n,m}$ is a unital morphism of algebras)
 we find that the Grothendieck bialgebra
is again dual to $\NCSym$, but this time the classes of simple modules
correspond to the monomial basis of $\NCSym$.

In section \ref{sec:frob} we give the proof of our first
theorem and show the basis canonically defined in section
\ref{sec:groth} corresponds to the simple modules of the $\kk \Pi_n$.
In light of the Frobenius characteristic of
section \ref{sec:groth}, the basis can
be interpreted as an analogue of the Schur functions
for $\NCSym$ and providing an answer to an open question
of \cite{RS}.

\section{$\NCSym$ and the basis $\{\X_A\}$}\label{sec:basis}

We recall the basic definition and structure of $\NCSym$.
Most of it can be found in \cite{BRRZ,BZ}.
A set partition $A$ of $m$ is a set of non-empty subsets
$A_1, A_2, \ldots, A_k \subseteq [m] = \{1,2, \ldots, m \}$
such  that $A_i \cap A_j = \emptyset$ for $i \neq j$ and
$A_1 \cup A_2 \cup \cdots \cup A_k = [ m ]$.
The subsets $A_i$ are called the parts of the set partition and
the number of non-empty parts  the length of $A$,
denoted by $\ell(A)$.  There is a natural mapping from set
partitions to integer partitions given by $\lambda(A)=( |A_1|,
|A_2|, \ldots, |A_k| )$, where we assume that the blocks of the
set partition have been listed in weakly decreasing order of size.
We shall use $\ell(\lambda)$ to refer to the length (the number
of parts) of the partition and $|\lambda|$ is the size of the
partition (the sum of the sizes of the parts), while
$n_i(\lambda)$ shall refer to the number of parts of the
partition of size $i$. We denote by $\Pi_m$
 the set of set partitions of $m$.  The number of set partitions is given
by the Bell numbers.  These can be defined
by the recurrence $B_0 = 1$ and
$B_n = \sum_{i=0}^{n-1} {n-1 \choose i} B_i$.

For a set $S = \{ s_1, s_2, \ldots, s_k \}$ of integers $s_i$ and
an integer $n$ we  use the notation $S+n$ to represent the set $\{
s_1+n, s_2+n, \ldots, s_k+n\}$. For $A \in \Pi_m$ and $B \in
\Pi_r$ set partitions with parts $A_i$, $1\leq i \leq \ell(A)$ and
$B_i$, $1 \leq i \leq \ell(B)$ respectively, we set $A \join B =
\{ A_1, A_2, \ldots, A_{\ell(A)}, B_1 + m, B_2+m, \ldots,
B_{\ell(B)}+m \}$, therefore $A \join B \in \Pi_{m+r}$ and this
operation is noncommutative in the sense that, in general, $A
\join B \neq B \join A$.

When writing examples of set partitions, whenever the
context allows it, we will use a more compact notation.
For example, $\{ \{1,3,5\}, \{2\}, \{4\} \}$
will be represented by $\{ 135\sep 2\sep 4\}$.  Although
there is no order on the parts of a set partition, we will
impose an implied order such that the parts are arranged
by increasing value of the smallest element in the subset.
This implied order will allow us to reference the $i^{th}$
parts of the set partition without ambiguity.

There is a natural lattice structure on the set partitions of
a given $n$.  We  define for $A, B \in \Pi_n$ that $A \leq B$
if for each $A_i \in A$ there is a $B_j \in B$ such that
$A_i \subseteq B_j$ (otherwise stated, that $A$
is finer than $B$).  The set of set partitions of $[n]$ with
this order forms a poset with rank function given by $n$ minus
the length of the set partition.
This poset has a unique minimal element
$\zero_n=\{1\sep 2\sep\,  \ldots\, \sep n\}$
and a unique  maximal element $\one_n=\{12\ldots n\}$.
The largest element smaller than both $A$ and $B$ is denoted
$A \wedge B = \{ A_i \cap B_j : 1\leq i \leq \ell(A), 1 \leq j \leq \ell(B) \}$
while the smallest element larger than $A$ and $B$ is denoted
$A \vee B$. The lattice $(\Pi_n, \wedge, \vee)$ is called the
 partition lattice.

\begin{exem}\textnormal{
Let $A = \{ 138\sep 24\sep 5\sep 67 \}$ and
$B = \{1\sep 238\sep 4567 \}$.  $A$ and $B$
are not comparable in the inclusion order on set partitions.  We
calculate that $A \wedge B = \{ 1\sep 2\sep 38\sep 4\sep 5\sep 67 \}$ and
$A \vee B = \{ 12345678 \}$.}
\end{exem}

When a collection of disjoint sets of positive integers
is not a set partition because
the union of the parts is not $[n]$ for some $n$, we may lower the values in
the sets so that they keep their relative values so that
the resulting collection
is a set partition (of an $m<n$).  This operation is
referred to as the `standardization' of a set of disjoint sets $A$ and the
resulting set partition is denoted $st(A)$.

Now for $A \in \Pi_m$ and $S \subseteq \{ 1, 2, \ldots, \ell(A)\}$ with
$S= \{ s_1, s_2, \ldots, s_k\}$, we define
$A_S = st(\{ A_{s_1}, A_{s_2}, \ldots, A_{s_k} \})$
which is a set partition of $|A_{s_1}| + |A_{s_2}| +
\ldots + |A_{s_k}|$.  By convention $A_{\{\}}$ is the empty set partition.

\begin{exem}\textnormal{
 If $A = \{ 1368\sep 2\sep 4\sep 579 \}$, then
 $A_{\{1,4\}} = \{ 1246\sep 357 \}$.}
\end{exem}

\medskip
For $n\ge 0$, consider a set $X_n$ of non-commuting variables
$x_1, x_2, \ldots, x_n$ and  the polynomial algebra
$\RR_{X_n} = {\kk}\langle x_1, x_2, \ldots, x_n \rangle$
in  these non-commuting variables.
There is a natural $S_n$ action on the basis elements defined
by $\sigma( x_{i_1} x_{i_2} \cdots x_{i_k} )
= x_{\sigma(i_1)} x_{\sigma(i_2)} \cdots x_{\sigma(i_k)}.$
Let $x_{i_1} x_{i_2} \cdots x_{i_m}$ be a monomial in the space
$\RR_{X_n}$.  We say that the type of this monomial is a set
partition $A \in \Pi_m$ with the property
that $i_{a} = i_{b}$ if and only if $a$ and $b$ are in the same
block of the set partition.  This
set partition is denoted as $\nabla(i_1, i_2, \ldots, i_m) = A$.
Notice that the length of $\nabla(i_1, i_2, \ldots, i_m)$ is
equal to the number of different values which appear in
$(i_1, i_2, \ldots, i_m)$.

The vector space $\NCSym^{(n)}$ is defined as the linear span of the
elements
$$\M_{A}[X_n] = \sum_{\nabla(i_1, i_2, \ldots, i_m)
= A} x_{i_1} x_{i_2} \cdots x_{i_m}$$
for $A \in \Pi_m$,
where the sum is over all sequences with $1 \leq i_j \leq n$.
For the empty set partition, we define by convention
$\M_{\{\}}[X_n] = 1$.
If $\ell(A)>n$ we must have that $\M_A[X_n] = 0$.
Since for any permutation $\sigma \in S_m$,
$\nabla(i_1, i_2, \ldots, i_m) =
\nabla(\sigma(i_1), \sigma(i_2), \ldots, \sigma(i_m))$, we have that
$\sigma \M_A[X_n] = \M_A[X_n]$. In fact, $\M_A[X_n]$ is the sum
 of all elements in the orbit of a monomial of type $A$
under the action of $S_n$. Therefore
 $\NCSym^{(n)}$ is the space of $S_n$-invariants in the noncommutative
 polynomial algebra $\RR_{X_n}$. For instance,
$ \M_{\{13\sep 2\}}[X_4]   =  x_1 x_2 x_1 + x_1 x_3 x_1 + x_1 x_4 x_1 + x_2 x_1 x_2
+ x_2 x_3 x_2 + x_2 x_4 x_2 + x_3 x_1 x_3 + x_3 x_2 x_3 + x_3 x_4 x_3
+ x_4 x_1 x_4 + x_4 x_2 x_4 + x_4 x_3 x_4.$

Now let $\NCSym$ be the inverse limit of $\NCSym^{(n)}$.
We have that $\NCSym = \bigoplus_{d\geq 0} \NCSym_d$ is a graded
algebra where $\NCSym_d $ is the linear span of
$\{ \M_A \}_{A \in \Pi_d}$.
Here we forget any reference to the variables $x_1,x_2,\ldots$
and think of elements in $\NCSym$ as noncommutative symmetric functions.
The degree of a basis element $\M_A$ is given by $|A|=d$
and the product map
$\mu : \NCSym_d \otimes \NCSym_m \longrightarrow \NCSym_{d+m}$
is defined on the basis elements $\M_A \otimes \M_B$ by
\begin{equation}  \label{eq:multm}
\mu( \M_A \otimes \M_B ) :=
\sum_{{{C \in \Pi_{d+m}}\atop
{C\,\,\wedge\,\,\one_d\join \one_m\, =\, A\join B}}} \M_C
\end{equation}
This is a lift of the multiplication in $\NCSym^{(n)}$ as $n\to\infty$.

The graded algebra $\NCSym$ is in fact a Hopf algebra
with the following comultiplication
$\Delta : \NCSym_d \longrightarrow \bigoplus_{k=0}^d
\NCSym_k \otimes \NCSym_{d-k}$
where
\begin{equation} \label{eq:comultm}
\Delta( \M_A ) = \sum_{S \subseteq [\ell(A)]} \M_{A_S} \otimes \M_{A_{S^c}}
\end{equation}
and $S^c = [\ell(A)] - S$. The counit is given by
$\epsilon\colon \NCSym\to {\mathbb Q}$
where $\epsilon(\M_{\{\}})=1$ and $\epsilon(\M_A)=0$
for all $A \in \Pi_n$ for $n>0$.
More details on this Hopf algebra structure are found
in \cite{BRRZ,BZ}.
The algebra $\NCSym$ was originally considered
by Wolf \cite{W} in extending the fundamental theorem
of symmetric functions to this algebra and later by
Bergman and Cohn \cite{BeCo}.
More recently Rosas and Sagan \cite{RS} considered this
space to define natural bases which are analogous to
bases of the (commutative) symmetric functions.  More
progress in understanding this space was made in
\cite {BRRZ,BZ} where it was considered as a Hopf algebra.
In the Hopf algebra $Sym$ of (commutative)
symmetric functions, the comultiplication corresponds to
the plethysm $f[X]\mapsto f[X+Y]$.
It was established in \cite{BRRZ} that the comultiplication in
$\NCSym$ corresponds to a noncommutative plethysm
$F[X]\mapsto F[X+Y]$.

The Hopf algebra $Sym$ has more structure. There is a second
comultiplication corresponding to the plethysm
$f[X]\mapsto f[XY]$ (see \cite{mac,Ze}). This second operation
is often refer to as the internal comultiplication or Kronecker comultiplication. We end this section describing for $\NCSym$
the analog of this internal comultiplication. This description
is also considered in \cite{AM}.

For the Hopf algebra $\NCSym$ we define a second (internal)
comultiplication
$\Delta^{\odot} : \NCSym_d \longrightarrow \NCSym_d \otimes \NCSym_d$
given by
\begin{equation} \label{eq:icomultm}
\Delta^{\odot}( \M_A ) = \sum_{B\wedge C = A} \M_{B} \otimes \M_{C}
\end{equation}
This operation corresponds to a noncommutative plethysm
$F[X]\mapsto F[XY]$. More precisely, we first assume that we have
two countable alphabet (totally ordered noncommutative variables)
$X=x_1, x_2,\ldots$ and $Y=y_1,y_2,\ldots$, and we let
$XY=\{x_iy_j\}$ be the alphabet totally ordered by the lexicographic
order on the couple $(i,j)$. The transformation $F[X]\mapsto F[XY]$
denotes the substitution given by the (countable) bijection
$X \leftrightarrow XY$. If we let the $x_i$'s commute with the
$y_j$'s then we have that $F[XY]$ can be expanded in the form
$F[XY]=\sum F_{1,i}[X]F_{2,i}[Y]$. We can then define the
operation $\Delta^{\odot}( F )=\sum  F_{1,i}\otimes F_{2,i}$.
In equation \eqref{eq:icomultm} we have given the result of this
when $F=\M_A$. Clearly this operation is a morphism for the
multiplication, thus $\NCSym$ with $\Delta^{\odot}$ and
the multiplication operation of equation \eqref{eq:multm} forms
a bialgebra. But it is not a Hopf algebra as it does not
have an antipode (the equation
$\mu\circ(\M_{\{\}}\otimes S)\circ\Delta^{\odot}(\M_{\{1\}})=\M_{\{\}}$
 has no solutions, where $S$ would have been the antipode).
 We are now in position to state our
first main theorem.

\begin{theo} \label{ThA} There is a basis
$\{\X_A : A \in \Pi_n, n\ge 0\}$ of $\NCSym$ such that
\begin{align*}
&\hbox{\rm (i)} \quad \X_A \X_B = \X_{A\join B}\hfill&\ &\ \\
&\hbox{\rm (ii)} \quad \Delta^{\odot}(\X_C) =
\sum_{A\vee B =C} \X_A\otimes \X_B .\hfill&\ &\
\end{align*}
\end{theo}

The proof of this theorem is technical and we differ it to
Section \ref{sec:frob}. We are convinced that the basis
$\{\X_A : A \in \Pi_n, n\ge 0\}$ is central in the study of
$\NCSym$ and should have many fascinating properties. We plan
to study this basis further in future work. For now, we prefer
to develop the representation theory that will motivate our result.

\section{Grothendieck bialgebras of the Semi-tower
$\bPiw =\bigoplus_{n\ge 0} \walg$.}\label{sec:groth}

In this section we consider the partition lattice algebras.
For a fixed $n$ consider the vector space $\walg$ formally
spanned by the set partitions of $n$. The multiplication is 
given by the operation $\wedge$ on set partitions and  with the unit 
$\one_n = \{ [n] \}$.
We remark that for all $d$, we have that  $\kk \Pi_d$
is isomorphic as a vector space to $\NCSym_d$ via the
pairing $A\leftrightarrow \M_A$. Moreover, it is
straightforward to check using equation \eqref{eq:icomultm}
that $\Delta^{\odot}$ is dual to $\wedge$ as operators.

It is well known that $\walg$ is a commutative semisimple
algebra (see \cite[Theorem 3.9.2]{stanley1}).  To see this, one
considers the algebra $\kk^{\Pi_n}=\{f\colon \Pi_n\to \kk\}$
which is clearly commutative and semisimple.
We then define the map
\begin{align*}
  \delta_{\ge}\colon \walg\   &\to  \ \kk^{\Pi_n} \\
                        A\  &\mapsto \ \delta_{A\ge},
\end{align*}
where $\delta_{A\ge}(B)=1$ if $A\ge B$ and $0$ otherwise.
Next check that $\delta_{A\wedge B\ge} =
\delta_{A\ge}\delta_{B\ge}$ which shows that $\delta_{\ge}$
is an isomorphism of algebras.

The primitive orthogonal idempotents of $\kk^{\Pi_n}$ are given by the
functions $\delta_{A=}$ defined by $\delta_{A=}(B)=1$ if
$A=B$ and $0$ otherwise. We have that
$\delta_{A\ge } = \sum_{B\le A} \delta_{=B}$. This implies,
using M\"obius inversion, that the primitive orthogonal idempotents of
$\walg$ are given by
\begin{equation}\label{minid}
  e_A=\sum_{B\le A} \mu(B,A) B,
\end{equation}
where $\mu$ is the M\"obius function of the partially ordered
set $\Pi_n$.
Since $\walg$ is commutative and semisimple,
we have that the simple $\walg$-modules of this algebra are
the one dimensional
spaces $V_A=\kk \Pi_n \wedge e_A$. Here the action is
given by the left multiplication
\begin{equation}\label{actionV}
 C\wedge e_A = \left\{ {e_A \quad\hbox{if $C\ge A$,} \atop 0\quad\hbox{otherwise.} } \right.
\end{equation}
This follows from the corresponding identity in $\kk^{\Pi_n}$ considering  $\delta_{\ge C}\delta_{=A}$.

We now let $G_0\walg$ denote the Grothendieck group of the category of
finite dimensional $\walg$-modules. This is the vector space
spanned  by the equivalence classes of simple $\walg$-modules
under isomorphisms.
We also consider $K_0\walg$  the Grothendieck
group of the category of  projective $\walg$-modules.
Since $\walg$ is semisimple, the space $G_0\walg$
and $K_0\walg$  are equal as vector spaces as they
are both linearly spanned by the elements $V_A$ for
$A \in \Pi_n$. We then
set  $K_0\bPiw=\bigoplus_{n\ge 0} K_0\walg$.

Given two finite $\walg$ modules $V$ and $W$,
we can form the $\walg$-module $V\otimes W$
with the diagonal action (it is an action
since a semigroup algebra is a bialgebra for the coproduct $A\to A\otimes A$).
  We denote this $\walg$-module by $V \odot W$
 (to avoid confusion with the tensor product of a $\walg$-module and a
 $(\kk\Pi_m,\wedge)$-module).
\begin{lemm}
Given two simple $\walg$-module $V_A$ and $V_B$,
\begin{equation}\label{tensorrep}
  V_A \odot V_B = V_{A\vee B}.
\end{equation}
\end{lemm}

\begin{proof}
Let $C\in\Pi_n$ act on $e_A\otimes e_B$.
>From equation \eqref{actionV} we get $C\wedge(e_A\otimes e_B)=
(C\wedge e_A) \otimes (C\wedge e_B)= e_A\otimes e_B$ if and only if $C\ge A$
and $C\ge B$, that is $C \ge A\vee B$.
If not, we get $C\wedge (e_A\otimes e_B)=0$. We conclude that the map
$e_A\otimes e_B\mapsto e_{A\vee B}$
is the desired isomorphism in equation \eqref{tensorrep}.
\end{proof}

We would like to define on $G_0\bPiw=\bigoplus_{n\ge 0} G_0\walg$
a graded multiplication and a graded comultiplication
corresponding to induction and restriction. For this we need a few
more tools.

\begin{lemm}\label{rhomult}
The linear map $\rho_{n,m}\colon \walg \otimes
(\kk \Pi_m,\wedge)\to (\kk \Pi_{n+m},\wedge)$
defined by $\rho_{n,m}(A\otimes B)=A|B$
is injective and multiplicative. Moreover,
$\rho_{k+n,m}\circ(\rho_{k,n}\otimes Id)=
\rho_{k,n+m}\circ(Id\otimes\rho_{n,m})$ for all $k,n$ and $m$.
\end{lemm}

\begin{proof}
Let $A=\{A_1,\ldots,A_r\}$, $B=\{B_1,\ldots,B_s\}$,
$C=\{C_1,\ldots,C_t\}$ and $D=\{D_1,\ldots,D_u\}$,
where $A,B\in\Pi_n$ and $C,D\in \Pi_m$. We remark that
for all $i,j$, we have $A_i\cap (D_j +n)=\emptyset$ and
$(C_i+n)\cap B_j=\emptyset$. Since $(C_i+n)\cap(D_j+n)=
(C_i\cap D_j)+n$, we have
\begin{align*}
(A|C)\wedge(B|D) &=
\big\{A_i\cap B_j \big\}_{1\le i\le r \atop 1\le j\le s}
\cup \big\{ (C_i+n)\cap(D_j+n) \big\}_{1\le i\le t \atop 1\le j\le u} \\
&= (A\wedge B)\,\big|\, (C\wedge D),
\end{align*}
and this shows that $\rho_{n,m}$ is multiplicative. The injectivity
of this map is clear from the fact that $\rho_{n,m}$ maps distinct
basis elements into distinct basis elements. The last identity of the lemma follows from the associativity of the operation ``$|$''
\end{proof}

We define a semi-tower $( \bigoplus_{n \geq 0} A_n, \{ \phi_{n,m}
\})$ to be a direct sum of algebras along with a family of
injective non-unital homomorphisms of algebras
 $\phi_{n,m} : A_n \otimes A_m \rightarrow A_{n+m}$.
A tower in the sense defined in the recent literature
\cite{BHT,Hiv,NCSF4} is a semi-tower with the additional
constraint that $\phi_{n,m}(\one_n, \one_m) = \one_{n+m}$ (i.e.
that $\phi_{n,m}$ is a unital embedding of  algebras).

Define the pair $\bPiw=\big(\bigoplus_{n\ge 0} \walg,\
\{\rho_{n,m}\}\big)$ which is a semi-tower of the algebras
$\walg$. We remark that $\bPiw$ is a graded algebra with the
multiplication $\rho_{n,m}(A,B)=A|B$ which is associative (but
non-commutative) and has a unit given by the emptyset partition
$\emptyset\in\Pi_0$. Moreover, each of the homogeneous components
$\walg$ of $\bf \Pi$ are themselves algebras with the
multiplication $\wedge$, and Lemma~\ref{rhomult} gives the
relationship between the two operations.

At this point we need to stress that $\rho_{n,m}$ is not a unital
embedding of algebras and hence $\bPiw$ is not a tower of
algebras. The algebra $\walg$ has a unit given by
$\one_n=\{12\ldots n\}$, but $\rho_{n,m}(\one_n \otimes \one_m)\ne
\one_{n+m}$. The tower of algebras considered in the recent
literature \cite{BHT,Hiv,NCSF4} all have the property that the
corresponding $\rho_{n,m}$ are (unital) embeddings of algebras.
This is the reason we call our construction a {\sl semi-tower}
rather than a tower.

The motivation for defining a tower of algebras is to allow one to
induce and restrict modules of these algebras and ultimately
to define on its Grothendieck ring a Hopf algebra structure. Here
the fact that we have only a semi-tower causes some problems in
defining restriction of modules.
Yet we can still define a weaker version of restriction in our
situation. Let $A$ and $B$ be two finite dimensional algebras
and let $\rho \colon A \to B$ be a multiplicative injective
linear map.  Given a finite $B$-module $M$, we define
$$\Res_\rho M = \{m\in M : \rho(\one_A)m=m\}\subseteq M.$$
In the case where $\rho$ is an embedding of algebras this
definition agrees with the traditional one. More on this general
theory will be found in \cite{huilan} but here we focus our
attention on $\bPiw$.

\begin{lemm}\label{restricV}
For $k\le n$ and a simple $\walg$-module $V_A\in G_0\walg$,
\begin{equation*}
\Res_{\rho_{k,n-k}} V_A = \begin{cases} V_A &
\hbox{if $A=B|C$ for $B\in\Pi_k$ and $C \in\Pi_{n-k}$ }\\
0&\hbox{\ otherwise.}
\end{cases}
\end{equation*}
\end{lemm}

\begin{proof}
We have that
$\rho_{n,m}(\one_k\otimes \one_{n-k})\wedge e_A=
(\one_k | \one_{n-k})\wedge e_A= e_A$ if
$\one_k | \one_{n-k} \ge A$, and $0$ otherwise.
The condition $\one_k | \one_{n-k} \ge A$
is equivalent to $A=B|C$ where $A|_{1,\ldots,k}=B$ and
$A|_{k+1,\ldots,n+k}=C$.
\end{proof}

We can now define a graded comultiplication on
$G_0\bPiw$ using our definition of restriction.
For $V\in G_0\walg$ let
\begin{equation}\label{eq:rescoprod}
\Delta(V) = \sum_{k=0}^n \Res_{\rho_{k,n-k}} V.
\end{equation}
It follows from Lemmas~\ref{rhomult} that this operation
is coassociative. 
For a simple module
$V_A\in G_0\walg$, Lemma~\ref{restricV} gives us
\begin{equation}\label{comultV}
  \Delta(V_A)=\sum_{A=B|C} V_B\otimes V_C.
\end{equation}
Now we extend $\odot$ to $G_0\bPiw$ by setting $V_A \odot V_B=0$
if $V_A$ and $V_B$ are not of the same degree.

\begin{prop}\label{bialgebra} $(G_0\bPiw,\odot,\Delta)$ is a bialgebra.
\end{prop}
\begin{proof} Let $A,B\in\Pi_n$. By equation \eqref{tensorrep},
 it is  sufficient to prove that
 $\Delta (V_{A\vee B})=\Delta (V_A)\odot \Delta (V_B)$.
Using equation \eqref{restricV} we can easily reduce the
 problem to the following assertion:
there are $C\in\Pi_k$, $D\in\Pi_{n-k}$ such that
$A\vee B = C|D$ if and only if
 there are  $E,E'\in \Pi_k$, $F,F'\in\Pi_{n-k}$ such that
  $A=E|F$ and $B=E'|F'$. This follows then from definitions.
\end{proof}

It is thus natural to give a notion to induced
modules dual to restriction in Lemma~\ref{restricV}.
\begin{lemm}\label{induceP}
For two simple modules
$V_A=\kk \Pi_n \wedge e_A \in G_0\walg$ and
$V_B=\kk \Pi_m \wedge e_B \in G_0(\kk \Pi_m, \wedge)$ we
define
$$\Ind_{n,m} V_A\otimes V_B=
\kk \Pi_{n+m} \otimes_{\kk \Pi_{n}\otimes\kk \Pi_{m}}(\kk \Pi_{n}\wedge e_A\otimes
\kk \Pi_{m}\wedge e_B),$$
where $\kk \Pi_{n}\otimes\kk \Pi_{m}$ is embedded into $\kk \Pi_{n+m}$
 via $\rho_{n,m}$.

There is a natural isomorphism such that
$$\Ind_{n,m} V_A\otimes V_B \cong
\kk \Pi_{n+m} \wedge \rho_{n,m}( e_A \otimes e_B).$$
We have
\begin{equation} \label{inductV}
  \Ind_{n,m} V_A\otimes V_B= V_{A|B}.
\end{equation}
\end{lemm}

\begin{proof}  Consider the following isomorphism which allows
us to naturally realize $\Ind_{n,m} V_A\otimes V_B$ as an element
of $G_0(\kk \Pi_{n+m}, \wedge)$.
\begin{align*}
\Ind_{n,m} V_A\otimes V_B &=
\kk \Pi_{n+m} \otimes_{\kk \Pi_{n}\otimes\kk \Pi_{m}}(\kk \Pi_{n}\wedge e_A\otimes
\kk \Pi_{m}\wedge e_B)\\
&= \kk \Pi_{n+m} \otimes_{\kk \Pi_{n}\otimes\kk \Pi_{m}}( e_A\otimes e_B)\\
&= \kk \Pi_{n+m} \wedge \rho_{n,m}(e_A \otimes e_B) \otimes_{\kk \Pi_{n}\otimes\kk \Pi_{m}}(\mathbf  1_n \otimes \mathbf 1_n)\\
&\cong \kk \Pi_{n+m} \wedge \rho_{n,m}(e_A \otimes e_B) .
\end{align*}
By linearity
$$\rho_{n,m}(e_A\otimes e_B)=
e_A| e_B=
\sum_{C\le A}\sum_{D\le B} \mu(C,A)\mu(D,B) C|D.$$
We now remark that $\{E: E\le A|B\}
=\{C|D: C|D\le A|B\}
= \{ C|D : C\le A ,\, D\le B\}$.
This is isomorphic to the cartesian product
$\{ C: C\le A\}\times \{ D : D\le B\}$.
Since M\"obius functions are multiplicative with
respect to cartesian product we have
$$\rho_{n,m}(e_A\otimes e_B)=
\sum_{E\le A|B} \mu(E,A|B) E = e_{A|B}.$$
\end{proof}

It is clear now that $\Ind_{n,m}$ defines on $G_0\bPiw$
a graded multiplication $V_A\otimes V_B\mapsto V_{A|B}$ that is dual
to the graded comultiplication of $\Delta$ defined on $G_0\bPiw$.
We also define an internal comultiplication on $G_0\bPiw$
dual to equation~\eqref{tensorrep} such that
$\Delta^\odot: G_0\walg  \rightarrow G_0\walg \otimes G_0\walg$.
For $C\in\Pi_n$ let
  \begin{equation}\label{intcomultP}
  \Delta^{\odot}(V_C) =\sum_{A\vee B=C} V_A\otimes V_B.
\end{equation}
The space $G_0\bPiw$ with its graded multiplication given by induction
and comultiplication $\Delta^\odot$ is a bialgebra,
 by duality and Proposition~\ref{bialgebra}.
 The
main theorem of this section is a direct corollary to Theorem~\ref{ThA}.

\begin{theo} \label{ThB}
The map $F\colon G_0\bPiw\to \NCSym$ defined by
$F(V_A)=\X_A$ is an isomorphism of bialgebras.
\end{theo}

\begin{proof}
$G_0\bPiw$ is endowed with a product given by \eqref{inductV} and
an inner coproduct given by \eqref{intcomultP}.  Since $\NCSym$ is
known to be a bialgebra satisfying the relations given in Theorem
\ref{ThA}, the map $F$ is an isomorphism.
\end{proof}

The map $F$ is called the Frobenius map for our semi-tower.
Along with Theorem~\ref{ThA}, it shows that the basis
$\X_A$ of $\NCSym$ are the only functions that correspond
to the classes of simple modules in $G_0\bPiw$. This
defines $\X_A$ uniquely (up to automorphism) and for this
reason we think of them as the Schur functions for the semi-tower
$\bPiw$ of the symmetric functions in non-commutative variables.

\begin{rema} \rm
In \cite{Brown}, Brown shows that $\walg$ is the
semisimple quotient of the Solomon-Tits algebra $ST_n$
(see \cite{Tits} ). It is easy to lift our semi-tower
structure from $\bf \Pi$ to ${\bf ST}=\bigoplus ST_n$
via Brown's support map. Then $G_0({\bf ST}, \wedge)$ and
$G_0\bPiw$ are isomorphic as bialgebras.
\end{rema}

In \cite{huilan}, the conditions under which a tower of algebras
${\bf A} = (\bigoplus_{n \geq 0} A_n, \rho_{n,m})$ defines a Hopf
algebra structure on the Grothendieck rings $G_0({\bf A})$ and
$K_0({\bf A})$ are considered. Under certain conditions one would
expect that the Grothendieck ring $G_0({\bf A})$ of finite modules
forms a Hopf algebra with the operations of induction and
restriction which is isomorphic to the graded dual of the
Grothendieck ring $K_0({\bf A})$ of projective modules.

For the tower of algebras we are
considering here, it is not the case that $G_0\bPiw$ forms
a Hopf algebra because the operations of induction and restriction
are not even compatible as a bialgebra structure.
We have shown that $G_0\bPiw$ and $K_0\bPiw$ are endowed naturally with
a product given by the notion of induction in equation
\eqref{inductV} and coproduct given by the notion of
restriction given in equation \eqref{eq:rescoprod}.  It is easily
checked that these operations do not form a Hopf algebra structure.

We have found however that here we have $G_0\bPiw$ endowed with the
operations of induction and restriction
is isomorphic by the graded dual to $K_0\bPiw$ also endowed with
the same induction and restriction operations. This is because
the operation of restriction on $G_0\bPiw$ is dual to
the operation of induction on $K_0\bPiw$ and induction
on $G_0\bPiw$ is dual as graded operations to restriction
on $K_0\bPiw$.  This remark can be observed through the
duality in equations \eqref{inductV} and \eqref{eq:rescoprod}.


\section{Grothendieck bialgebras of the Tower
$\bPiv =\bigoplus_{n\ge 0} \valg$.}\label{sec:grothII}

In this section we consider a second algebra related
to the partition lattice and show that there is an
additional connection with the algebra $\NCSym$.  Define
$\valg$ to be the commutative algebra
linearly spanned by the elements of
$\Pi_n$ and endowed with the product $\vee$.  This algebra
has as a unit the
minimal  element $\zero_{n} = \{1\sep2\sep\cdots\sep n\}$
of the poset $\Pi_n$ since $\zero_n \vee A = A$
for all $A \in \Pi_n$.

As we constructed the primitive orthogonal idempotents for $\walg$,
we proceed by defining in a similar manner
\begin{align*}
\delta_{\le} : \valg\ &\to\ \kk^{\Pi_n}\\
A\ &\mapsto\ \delta_{A\le}.
\end{align*}

It is straightforward to check that
$\delta_{A\le} \delta_{B\le} = \delta_{(A\vee B)\le}$
and hence $\delta_\le$ is an isomorphism of algebras.
This map can be used to recover the primitive orthogonal idempotents of $\valg$
since if $\delta_{A\le} = \sum_{B \ge A} \delta_{B=}$,
then $\delta_{A=} = \sum_{B \ge A} \mu(A,B) \delta_{B\le}$.
This can be summarized in the following proposition.

\begin{prop}
The primitive orthogonal idempotents of the algebra $\valg$ are
$$f_A = \sum_{B \ge A} \mu(A,B) B$$
with the property that
\begin{equation}\label{eq:actionfA}
C\vee f_A = \begin{cases}
f_A & \hbox{if }C\le A\\
0&\hbox{otherwise}
\end{cases}
\end{equation}
\end{prop}


It is also not difficult to check that the map
$\rho_{n,m}(A,B) = A|B$ is also multiplicative with
respect to the $\vee$ product in analogy with
Lemma \ref{rhomult}.  Therefore we define
the tower of algebras
$\bPiv = \big( \bigoplus_{n \geq 0} \valg, \{ \rho_{n,m} \}\big)$.
This time we find that $\rho_{n,m}$ is indeed an embedding of
algebras and $\bPiv$ a tower of algebras
(see remarks related to $\bPiw$)
since $\rho_{n,m}(\zero_n, \zero_m) =\zero_{n+m}$.

We now define $G_0\valg$ to be the ring of the category of finite dimensional
$\valg$-modules endowed with the tensor of modules as the product.
$G_0\valg$ is linearly spanned by the equivalence
classes of the simple modules under isomorphism.
  Also set $K_0\valg$ to be the
Grothendieck ring of the category of projective $\valg$-modules.
 Since $\valg$ is semi-simple we find that $G_0\valg \cong
K_0\valg$ and both are spanned by the simple modules
$W_A = \kk \Pi_n \vee f_A$.  Set $G_0\bPiv = \bigoplus_{n \ge 0} G_0\valg$
and $K_0\bPiv = \bigoplus_{n \ge 0} K_0\valg$.

Given two simple
modules $W_A, W_B \in G_0\valg$ we consider the
tensor product of modules $W_A \otimes W_B$ with the
diagonal action of $\valg$.  Denote this module as
$W_A \odot W_B$.  We find that
$C\vee ( f_A \otimes f_B) = (C\vee f_A) \otimes (C\vee f_B)$
which is equal to $f_A \otimes f_B$ if $C \le A$ and $C \le B$
(i.e. $C \le A \wedge B$)
and it is equal to $0$ otherwise.
We conclude from this discussion the following lemma.
\begin{lemm} For $W_A, W_B \in G_0\valg$,
\begin{equation}\label{eq:tensorpv}
W_A \odot W_B = W_{A \wedge B}
\end{equation}
is a simple module.
\end{lemm}

We have the following formula for
the restriction of $W_A$ to $\kk \Pi_k \otimes \kk \Pi_{n-k}$.

\begin{lemm}\label{restricW}
For $k\le n$ and a simple module $W_A\in G_0\valg$,
\begin{equation*}
\Res_{\rho_{k,n-k}} W_A =  W_B \otimes W_C
\end{equation*}
where $A \wedge (\one_k|\one_{n-k}) =B|C$
for $B\in\Pi_k$ and $C \in\Pi_{n-k}$.
\end{lemm}

\begin{proof}
First we check that
$\rho_{n,m}(\zero_k\otimes \zero_{n-k})f_A=
\zero_n \vee f_A= f_A$.
Now for $B' \in \Pi_k$ and $C' \in\Pi_{n-k}$, we have
that $\rho_{k,n-k}(B',C') \vee f_A = (B'|C') \vee f_A = f_A$ if
$(B'|C') \le A$ and $0$ otherwise.  If $A \wedge (\one_k|\one_{n-k})
 = (B|C)$
then $(B'|C') \vee f_A = f_A$ if and only if $B'\leq B$ and $C' \leq C$.
Therefore $W_A$ is isomorphic to $W_{B} \otimes W_{C}$ as
a $\kk \Pi_k \otimes \kk \Pi_{n-k}$ module.
\end{proof}

Define now a notion of induction for $K_0\bPiv$ (as the dual of $G_0\bPiv$).
For $A \in \Pi_n$ and $B \in \Pi_m$, the induced
$(\kk \Pi_{n+m}, \vee)$ module is
$$\Ind_{n,m} W_A \otimes W_B = \kk \Pi_{n+m}
\otimes_{\kk \Pi_n \otimes \kk \Pi_m} (W_A \otimes W_B)$$
where we consider $\kk \Pi_n \otimes \kk \Pi_m
\cong \rho_{n,m}(\kk \Pi_n \otimes \kk \Pi_m) \subseteq \kk \Pi_{n+m}$.

\begin{lemm} For $A \in \Pi_n$ and $B \in \Pi_m$ we have that
\begin{equation}\label{eq:inductv}
\Ind_{n,m} W_A \otimes W_B =
\sum_{C \wedge (\one_n|\one_m) = A|B} W_C.
\end{equation}
\end{lemm}

\begin{proof} By proceeding as in the proof of Lemma~\ref{induceP}, we get
$$
\Ind_{n,m} W_A \otimes W_B \cong \kk \Pi_{n+m} \wedge \rho_{n,m}(f_A \otimes f_B).
$$
Therefore we just have to prove
$$
\rho_{n,m}(f_A \otimes f_B) = \sum_{C \wedge (\one_n|\one_m) = A|B} f_C.
$$
Since M\"obius functions are multiplicative with
respect to cartesian product we have on the one hand
$$
\rho_{n,m}(f_A \otimes f_B) = \sum_{E|F\geq A|B} \mu(A|B,E|F) E|F.
$$
On the other hand
\begin{eqnarray*}
\sum_{C \wedge (\one_n|\one_m) = A|B} f_C &=& \sum_{C \wedge (\one_n|\one_m)=A|B}
                                                  \sum_{D\geq C} \mu(C,D)D\\
&=&\sum_{E|F\geq A|B} \sum_{C \wedge (\one_n|\one_m)=A|B}
\sum_{{\scriptstyle D \wedge (\one_n|\one_m)=E|F\atop\scriptstyle D\geq C}} \mu(C,D)D\\
&=& \sum_{E|F\geq A|B} \mu(A|B,E|F)E|F \\
&&+\sum_{E|F\geq A|B}
            \sum_{{\scriptstyle D \wedge (\one_n|\one_m)=E|F\atop\scriptstyle D\not=E|F}}
             \left(\sum_{{\scriptstyle C \wedge (\one_n|\one_m)=A|B\atop\scriptstyle C\leq D}}
              \mu(C,D)\right)D .
\end{eqnarray*}
The result follows then from the following equality
$$
\sum_{{\scriptstyle C \wedge (\one_n|\one_m)=A|B\atop\scriptstyle C\leq D}}  \mu(C,D)
=\sum_{A|B\leq C\leq D}  \mu(C,D)=0.
$$
\end{proof}

Induction and restriction define a graded product and coproduct on the space
of $G_0\bPiv = \bigoplus_{n \geq 0} G_0\valg$.
Define on the elements $N \in G_0(\kk \Pi_n, \vee)$ the operation
\begin{equation}\label{eq:comultresv}
\Delta( N ) = \sum_{k=0}^n \Res_{k,n-k} N
\end{equation}
and for $M \in G_0(\kk \Pi_m, \vee)$,
\begin{equation}\label{eq:multindv}
N \cdot M = \Ind_{n,m} N \otimes M.
\end{equation}
$G_0\bPiv$ with the operation $\Delta$ defines a coalgebra and
$G_0\bPiv$ with the product of \eqref{eq:multindv} defines an algebra structure.
It is easily checked that the product and coproduct on $G_0\bPiv$ are not compatible as a bialgebra structure.

It is interesting to note that $G_0\bPiv$ endowed with the tensor
product \eqref{eq:tensorpv} and the coproduct $\Delta$ from
equation \ref{eq:comultresv} does define a bialgebra.  To
highlight the relationship with $\NCSym$, we define an internal
coproduct $\Delta^\odot$ on $K_0\valg$ which is the natural dual
to equation \eqref{eq:tensorpv}.  That is we define a map
$\Delta^\odot: K_0\valg \rightarrow K_0\valg \otimes K_0\valg$
such that
\begin{equation}\label{eq:innercoprodv}
\Delta^\odot(W_A) = \sum_{B \wedge C = A} W_B \otimes W_C.
\end{equation}

We can now show with the following theorem that $K_0\bPiv$
is a bialgebra and the simple modules in $K_0\bPiv$
correspond to the $\M$-basis on $\NCSym$.
\begin{theo}  The ring $K_0\bPiv$ endowed with product $M \cdot N :=
\Ind_{n,m} M \otimes N$ and coproduct $\Delta^\odot$ of equation
\eqref{eq:innercoprodv} defines a bialgebra.  Moreover, the map
$F : K_0\bPiv \rightarrow \NCSym$ given by
$F( W_A ) = \M_A$ is an isomorphism of bialgebras.
\end{theo}

\begin{proof}
Recall that $\NCSym$ is a bialgebra linearly spanned by elements
$\M_A$ with the product defined by
$$\M_A \M_B = \sum_{C \wedge (\one_n | \one_m) = A|B} \M_C$$
and an inner coproduct defined by
$$\Delta^\odot(\M_A ) = \sum_{B \wedge C = A} \M_B \otimes \M_C.$$
Equations \eqref{eq:innercoprodv} and \eqref{eq:inductv} show that the map $F( W_A) = \M_A$
is an isomorphism of bialgebras.
\end{proof}



This construction that we have presented here in the last two
sections of defining an algebra from a lattice operation and
looking at the modules is something that can
be done in a more general setting and is a tool that can be
used to analyze other Hopf algebras.  This will be the subject
of future work.

\section{Existence of the $\X_A$ and Frobenius characteristic}
\label{sec:frob}

We now prove our Theorem~\ref{ThA}.
It is useful at this point to introduce an intermediate basis of $\NCSym$.
In \cite{RS}, an analogue of the power sum basis is given by
\begin{equation}\label{eq:defPA}
  \P_A= \sum_{B\ge A} \M_B.
\end{equation}
This basis has many nice properties.

\begin{lemm}\label{Pprop}
The set  $\{\P_A : A\in \Pi_n, n\ge 0\}$
forms a basis of $\NCSym$ such that
\begin{align*}
&\hbox{\rm (i)} \quad \P_A \P_B = \P_{A\join B}\hfill&\ &\ \\
&\hbox{\rm (ii)} \quad \Delta^{\odot}(\P_A) = \P_A\otimes \P_A .\hfill&\ &\
\end{align*}
\end{lemm}

\begin{proof}
By triangularity, it is clear that the set
forms a basis.
Now, for $A\in \Pi_n$ and $B\in \Pi_m$ we have
 $$\P_A \P_B = \sum_{C \geq A} \sum_{D \geq B} \M_C \M_D
= \sum_{C \geq A} \sum_{D \geq B} \quad
 \sum_{E \,\wedge \,\,\one_n\join \one_m = C|D} \M_E $$
Notice that we have that if
 $E \wedge \one_n\join \one_m= C|D$, then
$E \ge C|D \ge A|B$.
Conversely, if $E \ge A|B$, then we find unique $C$ and $D$
such that
$C|D= E\wedge\, \one_n\join \one_m \ge A|B\,\wedge\,
\one_n\join \one_m
= (A\wedge \one_n)|(B\wedge \one_m) = A|B$.
This implies that the sum is equal to
$$\P_A\P_B= \sum_{E \geq A|B} \M_E = \P_{A|B}.$$
For the second equality, we have
\begin{align*}
\Delta^\odot( \P_A ) &=  \sum_{B \geq A} \Delta^\odot ( \M_B )
 = \sum_{B \geq A} \sum_{C \wedge D = B} \M_C \otimes \M_D\\
&= \sum_{C \geq A} \sum_{B \geq A} \sum_{D: C \wedge D = B}  \M_C \otimes \M_D
 = \sum_{C \geq A} \sum_{D \geq A}  \M_C \otimes \M_D\\
&= \P_A \otimes \P_A
\end{align*}
\end{proof}

We finally define our basis. Let
\begin{equation}\label{def:xbasis}
\X_A= \sum_{B\le A} \mu(B,A) \P_B.
\end{equation}

By triangularity, the set $\{\X_A:A \Pi_n,\, n\ge 0\}$
is an integral basis of $\NCSym$.
We now see that this basis has the required properties.

  \begin{lemm}\label{Xprop}
\begin{align*}
&\hbox{\rm (i)} \quad \X_A \X_B = \X_{A\join B}\hfill&\ &\ \\
&\hbox{\rm (ii)} \quad \Delta^{\odot}(\X_C)
= \sum_{A\vee B =C} \X_A\otimes \X_B .\hfill&\ &\
\end{align*}
\end{lemm}

\begin{proof}
Using the same argument as in Lemma~\ref{induceP} we have
 \begin{align*}
\X_A\X_B&=\sum_{C\le A}\sum_{D\le B} \mu(C,A)\mu(D,B) \P_C\P_D\\
               &=\sum_{C\le A}\sum_{D\le B} \mu(C,A)\mu(D,B) \P_{C|D}
               = \sum_{E\le A|B} \mu(E,A|B) \P_E = \X_{A|B}.\\
\end{align*}
 This shows the first identity. For the second,
 the left hand side of (ii) is
 \begin{equation}\label{LHSi}
 \Delta^\odot(\X_C)=\sum_{E\le C}  \mu(E,C) \Delta^\odot(\P_E)=
   \sum_{E\le C}  \mu(E,C) \P_E\otimes \P_E,
 \end{equation}
 and the right hand side is
  $$
\sum_{A\vee B =C} \X_A\otimes \X_B= \sum_{A\vee B =C}\
\sum_{{E\le A \atop F\le B}} \mu(E,A) \mu(F,B) \  \P_E\otimes \P_F.
 $$
Let us isolate the coefficient of $\P_E\otimes \P_F$
in the sum above we get
\begin {align}
  T^C_{E,F} &=\sum_{{E\le A \le C\atop F\le B\le C}}
  \  \sum_{A\vee B =C}\mu(E,A) \mu(F,B)\label{TTT}\\
&= \sum_{F\le B\le C} \left( \sum_{{E\le A \le C\atop A\vee B =C}}
\mu(E,A) \right) \mu(F,B). \nonumber
\end{align}
By symmetry (interchanging the role of $E$ and $F$ if needed),
we may assume that $F\not< E$.
In \cite{stanley1}, Corollary 3.9.3 is dual to the following statement
\begin{equation*}\label{Stan393}
 \sum_{{ A \le \one_n\atop A\vee B =\one_n}}\mu(\zero_n,A) =
         \left\{ { \mu(\zero_n,\one_n) \quad
         \hbox{if ${ B=\zero_n\hfill}$,} \atop 0\quad\quad\qquad\ \hbox{otherwise.}\hfill } \right.
\end{equation*}
where, as usual, $\zero_n= \{1\sep 2\sep\,\ldots\,\sep n\}$.
This implies that the sum of in bracket in equation \eqref{TTT} is equal to
\begin{equation}\label{MMM}
 \sum_{{E\le A \le C\atop A\vee B =C}}\mu(E,A) =
         \left\{ { \mu(E,C) \quad\hbox{if ${ B=E\hfill}$,} \atop 0\quad\qquad\quad\hbox{otherwise.}\hfill } \right.
\end{equation}
This follows from the fact that $\mu$ is multiplicative
and in general the interval $[E,C]\subseteq \Pi_n$ is
isomorphic to a cartesian product of (smaller) partition
lattices (see Example 3.9.4 in \cite{stanley1}). If we substitute
this back in equation \eqref{TTT} we have two cases to consider.
When $F\ne E$, our assumption that $F\not< E$ prohibits the
possibility that $F\le B=E$. Thus we must always have
$B\ne E$ and in this case $T^C_{E,F}=0$.
When $F=E$, the only value of $B$ where equation \eqref{MMM}
does not vanish is when $B=E=F$ and we get
$T^C_{E,E}=\mu(E,C)\mu(E,E)=\mu(E,C)$.
If we compare this to equation \eqref{LHSi} we conclude
our proof of (ii).
\end{proof}

Notice that the character of the module (the trace of the matrix
module) $V_B$ from formula \eqref{actionV} is given by the formula
$\chi^{V_B}(A) = \delta_{B\le}(A)$ when $A \in \kk\Pi_n$ acts on
$V_B$. We observe that equation \eqref{def:xbasis} for $\X_A$
yields
$$\P_A = \sum_{B \le A} \X_B = \sum_{B} \chi^{V_B}(A) \X_B.$$
This means that the characters for the simple modules
for $\bPiw$ are encoded in the change of basis coefficients between
the $\P$ and $\X$ basis.

Similarly, the character of the module $W_B$ when acted on
by the element $A \in \kk \Pi_n$ are given by the formula
$\chi^{W_B}(A) = \delta_{B\ge}(A)$ from equation \eqref{eq:actionfA}.
Of course the defining
relation of the $\P$ basis from equation \eqref{eq:defPA} shows
that $$\P_A = \sum_{B} \chi^{W_B}(A) \M_B.$$
We observe in this formula that
the characters of the simple modules of
$\bPiv$ are encoded in the change of basis coefficients between
the $\P$ and $\M$ basis.

Both these formulas are in fairly close analogy with the formula
for the expansion for the power basis in the Schur basis in the
algebra of the symmetric functions. There the change of basis
coefficients are the characters of the simple modules of the
symmetric group.  This shows that the $\P$-basis which was defined
by Rosas and Sagan \cite{RS} does represent the analogue of the
power basis in the algebra of the symmetric functions and the $\X$
and the $\M$ bases encode in their coefficients the characters of
the modules that they represent.

\begin{rema}\textnormal{ One could also define a third algebra $(\kk \Pi_n, \oper)$ where
$A \oper B = \delta_{A=B} A$ and construct the simple modules
as we have done here for $\walg$ and $\valg$.
This same construction shows that the simple modules of this
algebra satisfy a tensor product, induction and restriction operations
which make the Grothendieck ring (of the category
 of the finite dimensional projective modules) for this algebra isomorphic again
to $\NCSym$ as a bialgebra where the simple modules behave as the
elements $\P_A \in \NCSym$ and $\P_A$ is defined in
\eqref{eq:defPA}.}
\end{rema}

\begin{rema}\textnormal{ Summary of bases in $\NCSym$.\\
 The $\mathbf m$ basis:
\begin{eqnarray*}
\M_A \M_B &=& \sum_{C\wedge(\one_{n}|\one_k) = A|B} \M_C\\
\Delta( \M_A ) &=& \sum_{S \subseteq [\ell(A)]} \M_{A_S} \otimes \M_{A_{S^c}}\\
\Delta^{\odot}( \M_A ) &=& \sum_{B \wedge C = A} \M_B \otimes \M_C
\end{eqnarray*}
The $\P$ basis:
\begin{eqnarray*}
\P_A \P_B &=& \P_{A|B}\\
\Delta( \P_A )& =& \sum_{S \subseteq [\ell(A)]} \P_{A_S} \otimes \P_{A_{S^c}}\\
\Delta^{\odot}( \P_A ) &=& \P_A \otimes \P_A
\end{eqnarray*}
The $\X$ basis:
\begin{eqnarray*}
\X_A \X_B &= &\X_{A|B}\\
\Delta^{\odot}( \X_A )& =& \sum_{B\vee C= A} \X_B \otimes \X_C
\end{eqnarray*}
It would be interesting to find a formula for $\Delta(\X_A)$. }
\end{rema}


\end{document}